\documentclass[preprint,12pt]{elsarticle}
\usepackage[utf8]{inputenc}

\usepackage{fullpage}

\usepackage{xcolor}
\usepackage{amsmath}
\usepackage{mathrsfs}
\usepackage{amssymb}
\usepackage{latexsym}
\usepackage{amsthm}
\usepackage[list=true]{subcaption}

\usepackage{tikz}

\newtheorem{theorem}{Theorem}
\newtheorem{lemma}[theorem]{Lemma}
\newtheorem{obs}[theorem]{Observation}
\newtheorem{corollary}[theorem]{Corollary}

\newtheorem{prop}[theorem]{Proposition}
\newtheorem{claim}{Claim}
\newtheorem{claim-a}{Claim}
\newtheorem{claim-b}{Claim}
\newtheorem{claim-c}{Claim}

\theoremstyle{definition}
\newtheorem{problem}[theorem]{Problem}

\newcommand{\brac}[1]{\left\lbrace #1 \right\rbrace}
\newcommand{\card}[1]{\left| #1 \right|}
\newcommand{\bs}{$\blacksquare$}

\usepackage{ifpdf}

\begin{document}

\begin{frontmatter}



\title{{Vertex arboricity of cographs\tnoteref{t1}}}

\tnotetext[t1]{This research was supported by }

\author[SFUM]{Sebasti\'an Gonz\'alez Hermosillo de la Maza}
\ead{sga89@sfu.ca}
\author[SFUCS]{Pavol Hell}
\ead{pavol@sfu.ca}
\author[CINVESTAV]{C\'esar Hern\'andez-Cruz\corref{cor1}}
\ead{cesar@cs.cinvestav.mx}
\author[SFUCS]{Seyyed Aliasghar Hosseini}
\ead{sahossei@sfu.ca}
\author[SFUCS]{Payam Valadkhan}
\ead{pvaladkh@sfu.ca}

\address[SFUM]{Department of Mathematics\\
Simon Fraser University\\
University Dr 8888\\
V5A 1S6, Burnaby, BC, Canada}
\address[SFUCS]{School of Computing Science\\
Simon Fraser University}
\address[CINVESTAV]{Departamento de Computaci\'on\\
Centro de Investigaci\'on y de Estudios Avanzados del IPN\\
}

\cortext[cor1]{Corresponding author}

\begin{abstract}
Arboricity is a graph parameter akin to chromatic number, in that
it seeks to partition the vertices into the smallest number of
sparse subgraphs. Where for the chromatic number we are partitioning
the vertices into independent sets, for the arboricity we want to
partition the vertices into cycle-free subsets (i.e., forests).
Arboricity is NP-hard in general, and our focus is on the arboricity
of cographs. For arboricity two, we obtain the complete list of
minimal cograph obstructions. These minimal obstructions do
generalize to higher arboricities; however, we no longer have a
complete  list, and in fact, the number of minimal cograph
obstructions grows exponentially with arboricity. We obtain bounds
on  their size and the height of their cotrees.

More generally, we consider the following common generalization of
colouring and partition into forests: given non-negative integers
$p$ and $q$, we ask if a given cograph $G$ admits a vertex partition
into $p$ forests and $q$ independent sets. We give a polynomial-time
dynamic programming algorithm for this problem. In fact, the
algorithm solves a more general problem which also includes several
other problems such as finding a maximum $q$-colourable subgraph,
maximum subgraph of arboricity-$p$, minimum vertex feedback set and
minimum $q$ of a $q$-colourable vertex feedback set.
\end{abstract}

\end{frontmatter}

\section{Introduction}
\label{sec:intro}

The {\em vertex-arboricity} of a graph $G$ is the minimum $p$ such that the
vertices of $G$ can be partitioned into  $p$ subsets each of which {\em
induces a forest}. We contrast this with the {\em chromatic number} of $G$,
which  is the minimum number $q$ such that the vertices of $G$ can be
partitioned into $q$ subsets each of which is {\em  independent}. Like the
chromatic number, determining the vertex arboricity of graphs is NP-hard
in general \cite{garey1979}. We focus our attention on the class of cographs,
where both problems are polynomial-time solvable. We define a common
generalization as follows. We say that a graph $G$ is $(p,q)$-{\em
partitionable} if the vertex set of $G$ can be partitioned into $p$ forests
and $q$ independent sets. This problem is NP-hard in general as well, as
long as $2p+q \geq 3$ (and is polynomial-time solvable otherwise). On the
other hand, it follows from \cite{damaschkeJGT14} that this problem also
has a polynomial-time algorithm for any $p, q$, when restricted to cographs.
Moreover, it follows from \cite{damaschkeJGT14}  that the number of minimal
obstructions for cograph $(p,q)$-partitionability is finite, for any $p, q$.
We investigate such minimal obstructions for $(p,0)$-partitionability of
cographs, i.e., for arboricity $p$. We give a complete  answer only for
arboricity $2$, and give some useful information for general $p$. We also
give a concrete dynamic programming algorithm to decide whether a cograph is
$(p,q)$-partitionable after the deletion of at most $r$ vertices.  This last
problem, allowing the deletion of vertices, is natural for the dynamic
programming algorithm, but it is an  interesting problem which can be
formulated as follows.

Let $p, q$ and $r$ be non-negative integers and let $G$ be a graph.
A $(p,q,r)$-{\em partition} of $G$ is a partition $(P, Q, R)$ of its
vertex set such that the subgraph induced on $P$ has vertex-arboricity
$p$, the subgraph induced on $Q$ is $q$-colourable, and $R$ has at most
$r$ vertices. We say that $G$ is $(p,q,r)$-{\em partitionable} if it
admits a $(p,q,r)$-partition, and we say that $G$ is a minimal
$(p,q,r)$-obstruction if it is not $(p,q,r)$-partitionable  but every
induced subgraph is. (When $r=0$, we simplify $(p,q,0)$ to $(p,q)$ in
all the notation.) Note that finding the  minimum $r$ such that $G$ is
$(0,q,r)$-partitionable is the well-known problem of finding the maximum
$q$-colourable subgraph; finding the minimum $r$ such that $G$ is
$(p,0,r)$-partitionable is the problem of finding the maximum subgraph
of arboricity $p$; finding the minimum $r$ such that $G$ is
$(1,0,r)$-partitionable is the minimum vertex feedback set problem;
finding the minimum $q$ such that $G$ is $(1,q,0)$-partitionable is the
problem of finding the smallest $q$ such that $G$ has a $q$-colourable
vertex feedback set.

If $G$ and $H$ are graphs, then we denote the disjoint union of $G$
and $H$ by $G + H$, and so, if $n$ is a positive integer, the disjoint
union of $n$ different copies of $G$ will be denoted by $nG$.   The
join of $G$ and $H$ will be denoted by $G \oplus H$.

A {\em cograph} is a graph than can be obtained recursively from
the following rules
\begin{itemize}
    \item $K_1$ is a cograph.

    \item If $G$ is a cograph, then $\overline{G}$ is a cograph.

    \item If $G$ and $H$ are cographs, then $G + H$ is a cograph.
\end{itemize}

There are many interesting characterizations of the family of
cographs \cite{corneilDAM3}, but there are two that are
particularly useful when dealing with minimal obstructions for a
hereditary property.   A graph is a cograph if and only if it is
$P_4$-free, if and only if the complement of any of its nontrivial
connected subgraphs is disconnected. Notice also that the
complement operation can be replaced by the join of two graphs
($G \oplus H$).

The rest of the article is organized as follows.   In Section
\ref{sec:200}, cographs that are $(2,0,0)$-partitionable are
characterized in terms of 7 minimal obstructions; some
families of cograph minimal obstructions for $(p,0,0)$-partitions
are studied.  In Section \ref{sec:1q0} we consider minimal
obstructions for $(1,q,0)$-partitions, and notice how these
partitions are related to the independent feedback vertex set
problem.   Although finite, the cograph minimal obstructions
for $(p,0,0)$-partition can be very large, both in size and in
number, Section \ref{sec:p00} is devoted to present lower and
upper bounds for these parameters, as well as an upper bound
on the height of the cotree of a minimal obstruction. A polynomial
algorithm to determine the arboricity of a cograph is presented
in Section \ref{sec:alg}. In Section \ref{sec:conc} we present
conclusions and related open problems.


\section{All Minimal Cograph Obstructions for Arboricity 2}
\label{sec:200}

Note that a graph has arboricity one if and only if it has no cycles.
Thus there are precisely two cograph minimal obstructions for arboricity
one, the cycles $C_3=K_3$ and $C_4=\overline{2K_2}$.

We now introduce a family of cographs $\mathcal{A}_2$ consisting of
$$\brac{K_5, \overline{3 K_3}, 2 K_3 \oplus \overline{K_2},
	2 \left(\overline{2 K_2}\right) \oplus \overline{K_3}, \overline{2 K_2}
	\oplus (K_1 + K_2), \left(\overline{2 K_2} + K_3\right) \oplus
	\overline{K_2},  \overline{3 K_2 + K_1}}.$$

\begin{figure}[ht!]
	\centering
	\subcaptionbox{$K_5$}{\includegraphics{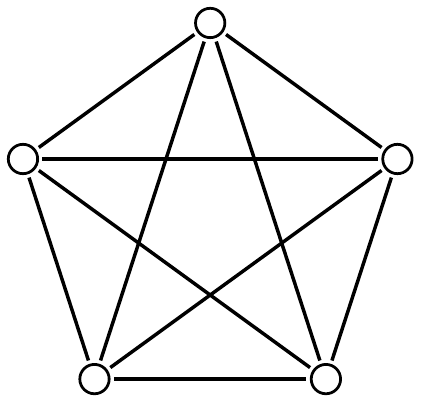}}%
	\hfill
	\centering
	\subcaptionbox{$\overline{3 K_3}$}{\includegraphics{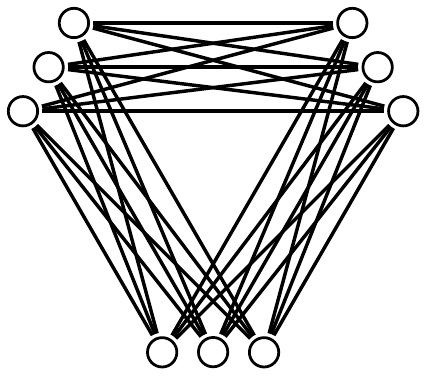}}%
	\hfill
	\centering
	\subcaptionbox{$\overline{2 K_2}
		\oplus (K_1 + K_2)$}{\includegraphics{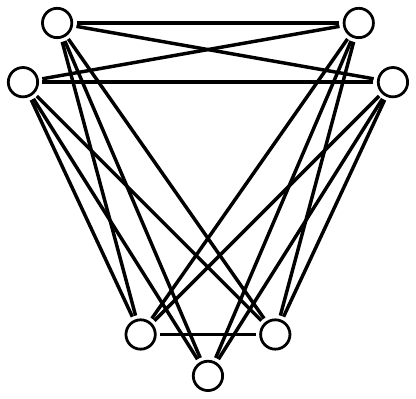}}%
	\hfill
	\centering
	\subcaptionbox{$2 \left(\overline{2 K_2}\right) \oplus \overline{K_3}$}{\includegraphics{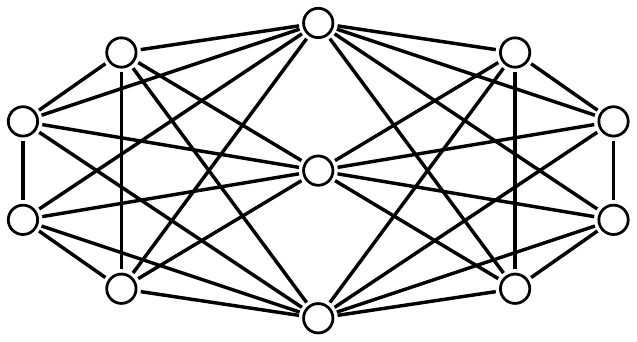}}%
	\hfill
	\vspace{5mm}
	\centering
	\subcaptionbox{$2 K_3 \oplus \overline{K_2}$}{\includegraphics{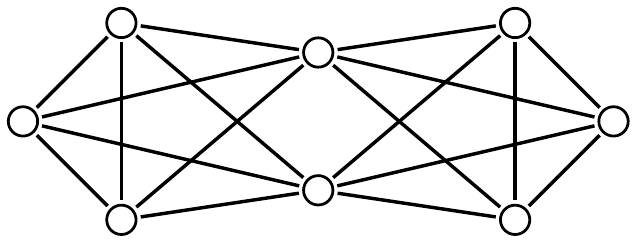}}
	\hfill
	\centering
	\subcaptionbox{$\overline{3 K_2 + K_1}$}{\includegraphics{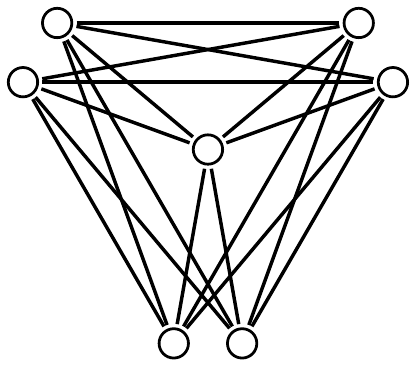}}
\hspace{10mm}
		\centering
	\subcaptionbox{$\left(\overline{2 K_2} + K_3\right) \oplus
		\overline{K_2}$}{\includegraphics{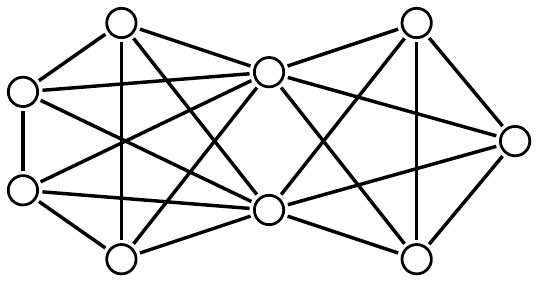}}
	\caption{The family $\mathcal{A}_2$.}\label{A2}
\end{figure}

These graphs are depicted in Figure \ref{A2}.

%
%
%
%
%
%
%
%
%

\begin{lemma}\label{xxx}
	Each graph in $\mathcal{A}_2$ is a cograph minimal obstruction
	to arboricity $2$.
\end{lemma}

\begin{proof}

	It is clear from the descriptions that each graph in $\mathcal{A}_2$
	is a cograph. We claim that each of them is not partitionable into
	two forests, but whenever a  vertex is deleted, it becomes so
	partitionable.

	We prove the first, fourth and last cases ($K_5$, $2 (\overline{2
	K_2}) \oplus \overline{K_{3}}$, $\overline{3 K_2 + K_1}$); the rest
	of the cases can be handled similarly. Consider first $G = K_5$. It
	is clear that in a complete graph each (acyclic) colour class has at
	most $2$ vertices, therefore $G$ is a $(2,0,0)$-obstruction. To check
	the minimality, remove any vertex of $G$, the remaining graph is a
	$K_{4}$ which is easily $2$-colourable. Therefore, $G$ is actually
	a minimal $(2,0,0)$-obstruction.

	Let us assume that $G = 2 (\overline{2 K_2}) \oplus
	\overline{K_{3}}$. First notice that to colour $\overline{K_{3}}$
	by using only $2$ colours, at least two vertices receive the same
	colour, and then we can use that colour on at most one other vertex
	outside $\overline{K_{3}}$. Since we cannot use this colour anywhere
	else, without loss of generality we can assume that all vertices of
	$\overline{K_{3}}$ have the same colour and we are using this colour
	on one other vertex as well. So we have coloured (at most) $4$
	vertices using one colour. There are two disjoint copies of
	$\overline{2 K_2}$ (minus one vertex) still uncoloured. On each of
	the $\overline{2 K_2}$ copies we can use one colour for at most $3$
	vertices. Therefore, at most $6$ vertices can be coloured, and thus
	at most $4 + 6 = 10$ vertices can be coloured with $2$ colours.
	Since $G$ has $11$ vertices, we conclude that $G$ is a
	$(2,0,0)$-obstruction.
	To verify minimality, first consider the case where $v\in \overline{K_{3}}$.
	Then we can colour one vertex of $\overline{K_{3}}$ along with one copy of
	$\overline{K_2}$ in each of the two $\overline{2 K_2}$ parts to colour $G
	- v$. Now consider the case where $v$ belongs to $2 (\overline{2 K_2})$.
	Let $v'$ be the duplicate of $v$ in the other copy of $\overline{2 K_2}$
	and colour $v'$ along with all vertices of $\overline{K_{3}}$ using one
	colour.   The remaining vertices induce an acyclic graph and can be coloured
	with one colour. This shows that $G - v$ is $2$ (acyclic) colourable.
	Therefore $G$ is a minimal $(2,0,0)$-obstruction.

	For our final case, let $G = \overline{3 K_2 + K_1}$. Notice that
	$G$ has $7$ vertices and each (acyclic) colour class has size at most
	3 (two vertices of a $\overline{K_2}$ and another vertex). Therefore
	$G$ is a $(2,0,0)$-obstruction. By removing any vertex of $G$ we will
	have $6$ vertices and at least $2$ copies of $\overline{K_2}$.
	Considering each of these $\overline{K_2}$'s along with some other
	vertex (maybe in another $\overline{K_2}$) will give us an acyclic
	colouring where each colour class has size exactly 3. Thus, $G$ is
	actually a minimal $(2,0,0)$-obstruction.
\end{proof}

\begin{theorem}\label{yyy}
A cograph has vertex arboricity $2$ if and only if it is
$\mathcal{A}_2$-free.
\end{theorem}

\begin{proof}
Let $G$ be a cograph.   If the vertex arboricity of $G$
is at most $2$, then it is clearly $\mathcal{A}_2$-free.
So, suppose that $G$ is $\mathcal{A}_2$-free.   We
may assume without loss of generality that $G$ is
connected.

Since $G$ is a connected cograph, there exist cographs
$G_1$ and $G_2$ such that $G = G_1 \oplus G_2$.   If
$G_1$ and $G_2$ are forests, we are done.   So, at least
one of them must contain an induced cycle. Without loss of
generality suppose that it is $G_1$.  Since $H$ is a cograph,
this cycle should be a triangle or a $4$-cycle.  Suppose that
$G_1$ is triangle free, then it must contain a $4$-cycle.
Since $G$ is $K_5$-free, then $G_2$ is triangle-free, and,
since $G$ is $(\overline{3 K_2 + K_1})$-free, then $G_2$ is
$P_3$-free, i.e., $G_2$ is a disjoint union of $K_1$'s and
$K_2$'s; but $G$ is $(C_4 \oplus (K_1 + K_2))$-free, then $G_2$
is $(K_1 + K_2)$-free, and thus, it is either a $K_2$, or
an empty graph.   As a triangle-free cograph, $G_1$ is
bipartite, with bipartition $(X,Y)$.   If $|V(G_2)| \le 2$,
colour red one vertex in $G_2$, together with all the vertices
in $X$, and colour blue the other vertex in $G_2$ (possibly
none), together with the vertices in $Y$, this colouring of
$G$ realizes the vertex arboricity $2$.   If $|V(G_2)| \ge 3$,
as $G$ is $(2 C_4 \oplus 3 K_1)$-free, we have that every
component of $G_1$, different from the one containing the
induced $4$-cycle, is a star.   The component of $G_1$
containing the induced $4$-cycle is a bipartite connected
cograph, and thus, it is a complete bipartite graph;
moreover, since $G$ is $\overline{3 K_3}$-free, one of the
two parts of this component has less than three vertices.
colour red one of the vertices in this small part, together with
all the vertices in $G_2$, and colour blue all the remaining
vertices of $G$.   Clearly, the red vertices induce a star, and
the blue vertices induce a disjoint union of stars.

Now, suppose that $G_1$ contains an induced triangle.
Using again that $G$ is $K_5$-free, we conclude that $G_2$
is an empty graph.   Let the set $\{ v_1, v_2, v_3 \}$ induce
a triangle in $G_1$, and let $B$ be the component of $G_1$
containing it.   Now, $B$ is a connected cograph, and there
are cographs $B_1$ and $B_2$ such that $B = B_1 \oplus
B_2$. Recall that $G_1$ is $K_4$-free, and thus, both
$B_1$ and $B_2$ are triangle-free, so, we assume without
loss of generality that $v_1, v_2 \in V(B_1)$, and $v_3 \in
V(B_2)$; moreover, $B_2$ must be an independent set.

We will consider two cases; suppose first that $G_2$ has
at least two vertices.   Then, since $G$ is $((K_3 + C_4)
\oplus 2 K_1)$ and $(2 K_3 \oplus 2 K_1)$-free, we have
that $G_1$ has precisely one component, namely $B$,
which is not acyclic.    Again, we have two cases. First,
suppose that $B_2$ has at least two vertices.   From the
fact that $G$ is $(C_4 \oplus (K_1 + K_2))$-free, we
obtain that $B_1$ is connected, and, since $G$ is
$(\overline{ 3 K_2 + K_1 })$-free, then $B_1$ is a path
on two vertices, $x$ and $y$.   So, we can colour $x$,
together will all the vertices in $G_2$ red, and the rest
of the vertices in $G$ blue; it is not hard to verify that
each colour class induces a forest.   So, we may
now suppose that $B_2$ has only one vertex.   Since
$B_1$ is triangle-free, it is bipartite with bipartition $(X,
Y)$.   Again, as $G$ is $(\overline{ 3 K_2 + K_1})$-free,
we have that $B_1$ is acyclic, and we can colour the only
vertex in $B_2$ together with all the vertices in $G_2$ red,
and the rest of the vertices in $G$ blue; again, each
colour class induces a forest.

As a second case, suppose that $G_2$ has only one vertex
$v$.   Since $G$ is $\{ K_5, \overline{3 K_2 + K_1} \}$-free,
and $v$ is a universal vertex in $G$, we have that $G_1$ is
$\{ K_4, \overline{3 K_2} \}$-free.   It follows from Theorem
\ref{ifvs} with $q=1$ that $G_1$ contains an independent feedback
vertex set, $S$.   Thus, colouring $v$ together with the
vertices in $S$ red, and the rest of the vertices of $G$
blue, clearly yields acyclic colour classes.
\end{proof}

The family $\mathcal{A}_2$ has a natural generalization for higher arboricity.
Let $p$ be an integer, $p \ge 2$, and denote by $\mathcal{A}_p$ the following
family of cographs.

\begin{itemize}
	\item $K_{2p+1}$

	\item $\overline{(p+1) K_{p+1}}$

	\item $2 K_{2p-1} \oplus \overline{K_2}$

	\item $2 \overline{p K_p} \oplus \overline{K_{p+1}}$

	\item $\overline{p K_2} \oplus (K_1 + K_p)$

	\item $(\overline{p K_{p}} + K_{2p-1}) \oplus
		\overline{K_2}$

	\item $\brac{ \overline{(p+1+i)K_2 + (p-1-2i)K_1}}_{0
		\leq i \leq \lfloor\frac{p-1}{2}\rfloor}$
\end{itemize}

\begin{lemma}
Let $p$ be an integer, $p \ge 2$.   Each graph in the family
$\mathcal{A}_p$ is a cograph minimal obstruction for arboricity $p$.
\end{lemma}

\begin{proof}
The proof is analogous to that of Lemma \ref{xxx}.
\end{proof}

There is however, no analogue to Theorem \ref{yyy}. In fact,
in Section \ref{sec:p00} we prove that the number of cograph
minimal obstructions for arboricity $p$ grows exponentially
with $p$.

\section{Minimal Cograph Obstructions for $q$-Colourable Vertex
Feedback Set}
\label{sec:1q0}

By analogy with the independent vertex feedback set problem, we
say that a cograph $G$ has a $q$-{\em colourable vertex feedback
set} if it has a $(1,q)$-partition. It turns out there are exactly
two minimal cograph obstructions for $(1,q)$-partition, namely, the
complete graph $K_{q+3}$, and the complete $(q+2)$-partite graph with
two vertices in each part, $\overline{(q+2)K_2}$.

\begin{theorem} \label{ifvs}
Let $q$ be a non-negative integer.   A cograph $G$ has a
$q$-colourable vertex feedback set if and only if it is
$\brac{K_{q+3}, \overline{(q+2)K_2}}$-free.
\end{theorem}

\begin{proof}
Clearly, a $(1,q,0)$-partitionable cograph is $\brac{K_{q+3},
\overline{(q+2)K_2}}$-free. We prove the converse by
induction on $q$. The base case $q=0$ follows from the
simple fact that a cograph is a forest is and only if it is
$\brac{K_3, \overline{2K_2}}$-free. Suppose that
the claim holds for all $\ell < q$, and let $G$ be a
$\brac{K_{q+3}, \overline{(q+2)K_2}}$-free cograph. Without
loss of generality, we may assume $G$ is connected. The
fact that $G$ is $K_{q+3}$-free implies $\chi(G) \leq q+2$,
and the claim holds by taking an independent set as the
forest if $\chi(G) \leq q+1$, so we may assume $\chi(G) =
q+2$.

Since $G$ is a connected cograph, there exists a family
of cographs $\brac{G_i}_{i=1}^s$ such that $G =
\bigoplus_{i=1}^s G_i$. Notice that in any $(q+2)$-colouring
of $G$, each colour class is contained in $V(G_i)$ for
some $i \in \brac{1, 2, \dots, s}$. If there is a $(q+2)$-colouring
of $G$ with a colour class with a single vertex $v$, then
$v$ together with any other colour class induces a forest.
By taking this forest and the remaining $q$ colour classes
we obtain a $(1,q,0)$-partition of $G$, so we may assume
that every colour class has at least two vertices.

Since $\sum_{i=1}^s \chi(G_i) = q+2$, if $G_i$ has an
induced $\overline{\chi(G_i)K_2}$ for every $i \in \brac{1,
\dots, s}$, then $G$ has an induced copy of
$\overline{(q+2)K_2}$, thus $G_\ell$ is $\left( \overline{
\chi(G_j)K_2}\right)$-free for some $\ell \in \brac{1, \dots, s}$.
By induction hypothesis, $G_\ell$ has a $(1,\chi(G_\ell)-2,
0)$-partition (notice that if $\chi(G_\ell) = 1$, the fact that each
colour class has at least two vertices would imply the
existence of an induced $\overline{K_2}$, and so we must
have that $\chi(G_\ell) \geq 2$). Since $G' = G - V(G_\ell)$ has
chromatic number $q+2-\chi(G_\ell)$, a proper colouring of
$G'$, together with the $(1, \chi(G_\ell)-2,0)$-partition of
$G_\ell$ gives us a $(1,q,0)$-partition of $G$.
\end{proof}

We can use the Theorem to derive a min-max relationship.
For the purposes of its statement, we shall call $K_s$ a
{\em thin $s$-clique} and $\overline{sK_2}$ a {\em thick $s$-clique}.
The {\em strength} of a thin $s$-clique is defined to be $s$, and the
strength of a thick $s$-clique is defined as $s+1$. We let $s(G)$
denote the maximum strength of a (thin or thick) clique in $G$. We
also let $q(G)$ denote the minimum number of colours $q$ such that
$G$ admits a $q$-colourable vertex feedback set.

\begin{corollary}
If $G$ is a cograph, then $q(G) = s(G) - 2$.
\end{corollary}

We note that the maximum strength of a clique in a cograph can be
computed by a cotree bottom-up procedure analogous to the well-known
algorithm for computing the size of a maximum complete subgraph of a
cograph \cite{corneilDAM3}.





\section{Bounds on Minimal Cograph Obstructions for Arboricity $p$}
\label{sec:p00}

As we mentioned at the end of Section \ref{sec:200}, we do
not have a complete description for all cograph minimal
obstructions for arboricity $p$, $p > 2$.   In this
section we illustrate the fact that there are exponentially
many.   We will construct this obstructions as joins of
star forests.

\begin{prop}
Let $p$ be an integer, $p \ge 2$.  There are at least
$\frac{e^{2 \cdot \sqrt{p}}}{14}$ minimal obstructions for
arboricity $p$.
\end{prop}

\begin{proof}
Let us first observe that if we add some edges to a minimal
obstruction, the resulting graph is still an obstruction but
might not be minimal. Consider the complete multipartite
graph $\overline{(p+1) K_{p+1}}$ (which is a minimal
obstruction for arboricity $p$), and add edges to  $i$ of
the parts to make them non-empty forests.

Let $\mathcal{F}_p$ be the set of all non-empty forests on
$p$ vertices which are cographs. Notice that a tree cograph
on a fixed number of vertices is unique (it is a star).
Therefore, $|\mathcal{F}_p|=\pi(p)-1$, where $\pi (p)$ is
the is the partition number of $p$ (the number of possible
partitions of $p$) . In \cite{marotiEJCNT3} the following
lower bound is proved for $\pi$ \[
    \frac{e^{2 \cdot \sqrt{p}}}{14} < \pi(p).
\]

Let $i$ be an integer, $0 \le i \le p$, and define the graph
$\mathcal{O}_i(f_1,\ldots,f_i)$ as $\mathcal{O}_i(f_1, \ldots,
f_i) = \overline{(p+1-i) K_{p+1-i}} \oplus f_1 \oplus f_2
\ldots \oplus f_{i}$, where $f_j \in \mathcal{F}_{p+2-i}$
for $j \in \{ 1, \ldots, i \}$. Notice that when $i=0$,
$\mathcal{O}_i$ does not receive arguments, and we obtain
$\overline{(p+1) K_{p+1}}$, and when $i=p$ we have that
$f_j$ is isomorphic to $K_2$ for every $j \in \{ 1, \dots,
p \}$, and thus, the only possible graph that can be obtained
is $K_{2p+1}$. We denote the set of all graphs $\mathcal{O}_i(f_1,
\ldots, f_i)$ for all selections of $f_j \in \mathcal{F}_{p+2-i}$
by $\mathcal{O}_i$.

\begin{claim}
    For every $i \in \{ 0, \ldots, p \}$, each $G\in \mathcal{O}_i$
    is not $p$-colourable.
\end{claim}
Let $G$ be a member of $\mathcal{O}_i$. The number of vertices
of $G$ is $(p+1-i)^2+(i)(p+2-i)=p^2+2p-ip+1=p(p+2-i)+1$. On the
other hand, notice that we can use each colour class in at most
two parts of $G$ (otherwise we will get a monochromatic cycle).
Also if we want to use a colour in two parts, in one of them we
are using it at most once (or we will get a monochromatic $C_4$).
Since each $f_j$ is non-empty, if we colour it with just one
colour, then we cannot use that colour for any other vertices
(or we will get a monochromatic triangle). Therefore, in each
colour class we can have at most $p+2-i$ vertices. Hence, using
$p$ colours we can colour at most $p(p+2-i)$ vertices of $G$.
Therefore $G$ is not $p$-colourable for $i \in \{ 0, \ldots, p
\}$. \bs

\begin{claim}
    For every $i \in \{ 0, \ldots, p \}$, every $G \in
    \mathcal{O}_i$ is a minimal obstruction for
    arboricity $p$.
\end{claim}
We have proved that $G$ is an obstruction, now we need to show
that it is minimal. If we remove a vertex from $f_j$, then we
can use each colour for one vertex of $f_j$ and one of the parts
of $\overline{(p+1-i) K_{p+1-i}}$, which is an independent set
(using $p+1-i$ colours). Also we can use one colour for each of
the remaining $f_l$ (using $i-1$ colours).

If we remove a vertex from a part of $\overline{(p+1-i)
K_{p+1-i}}$, then we can colour the remaining graph by using
only $p$ colours; use one colour for each vertex of this part
together with all the vertices in another of the parts, i.e.,
use each vertex of this part as the center of a star having
all the vertices in some other part as leaves (using $p-i$
colours) and use one colour for each $f_j$ (using $i$ colours).
\bs

Now we need to calculate $|\mathcal{O}_i|$. For a fixed value
of $i$, we have to consider the effect of permutation (switching
$f_i$ and $f_j$) by dividing each term by $i!$.  It is easy to
see that for two different sets of forests we will get different
obstructions. So we have \[
  |\mathcal{O}_i| = \frac{(\pi(p+2-i)-1)^{i}}{i!}.
\]
Notice that if $i \ne j$, then members of $\mathcal{O}_i$ and
$\mathcal{O}_j$ have different number of vertices and therefore
they are different. So the total number of different minimal
obstruction that we will get from this structure is \[
|\mathcal{O}_i| = \sum_{i=0}^{p}\frac{(\pi(p+2-i)-1)^i}{i!} \ge
\sum_{i=0}^{p}\frac{e^{2i\sqrt{p+2-i}}}{14^{i}i!} >
\frac{e^{2\sqrt{p}}}{14}.
\]
\end{proof}

%
%


Next we focus on upperbound for minimal cograph obstructions
to arboricity $p$.   First, we consider the cotree height.

\begin{theorem}
Let $p$ be an integer, $p \ge 2$.   If $G$ is a minimal cograph
obstruction for arboricity $p$ with cotree $T$, then the height
of $T$ is at most $4p+1$.
\end{theorem}

\begin{proof}
We may assume that  $G$ is connected and therefore the root
vertex of the cotree is a join vertex, $J_0$. Notice that to
have a unique co-tree all children of a join vertex should be
union or single vertices and all children of union vertices
should be join or single vertices. Recall that $K_{2p+1}$,
whose cotree has height two, is a minimal obstruction and
therefore no other minimal obstruction can contain it. For
simplicity, in the following  we will use $\omega(X)$ instead
of $\omega(G[X])$ to denote the clique number of the subgraph
of $G$ induced by the vertex set $X$.

Let $J$ be a join vertex of the co-tree with degree $d$ whose
children are $U_1, U_2, \ldots, U_d$. It is easy to see that
$\omega(J)  = \sum_{i=1}^d \omega \left(U_i\right)$. In
particular, $d \geq 2$ implies that $\omega(J) \geq \omega(U_i)
+ 1$ for any $U_i$ that is a child of $J$. This implies that
any path from $J_0$ to a leaf of $T$ contains at most $2p$ join
vertices, hence the height of the co-tree is at most $4p+1$.
\end{proof}

\begin{corollary}
\label{joinVertex}
Let $p$ be an integer, $p \ge 2$.   Let $G$ be a minimal cograph
obstruction for arboricity $p$ and let $T$ be its cotree.

If $G \neq K_{2p+1}$, then every join vertex in $T$ has at most
$2p$ children.
\end{corollary}

%

\begin{theorem}\label{exp}
Let $G_1$ and $G_2$ be minimal cograph obstructions for
$p$-vertex arboricity such that $\rho (G_i) = \chi(G_i)
= p+1$ for $i \in \brac{1,2}$. Let $T(G)$ denote the
height of the cotree of a cograph $G$. If $S$ is an
independent set of size  $p+2$, then the cograph $H =
(G_1 + G_2) \oplus S$ satisfies:

\begin{enumerate}

\item[(a)] $\rho(H) = \chi(H) = p+2$.

\item[(b)] $H$ is a cograph minimal obstruction for
$(p+1)$-arboricity.

\item[(c)] $T(H) = \max \brac{T(G_1), T(G_2)} +1$
\end{enumerate}

\end{theorem}

\begin{proof}
It is easy to see that $\chi(H)=p + 2$, and $\rho(H)
\leq p+2$. If $\rho(H) \leq p+1$, then there exists
a partition $\mathcal{F}$ of $V(H)$ into $p+1$ induced
forests, and hence, at least one forest $F$ in $\mathcal{F}$
contains two distinct vertices of $S$. This implies that
there exists $i \in \brac{1,2}$ such that $V(F) \cap V(G_i)
= \varnothing$, and so the restriction of $\mathcal{F}$ to
$V(G_i)$ is a partition of $V(G_i)$ into $\rho(G_i)-1 = p$
forests, which is a contradiction, and so $a$ holds.

To show $b$, let $v \in V(H)$. If $v \in S$, let $S'$ be $S'
= \brac{v_1, \dots, v_{p+1}} = S - \brac{v}$, and take
a $(p+1)$-colouring, $f_i \colon V(G_i) \to S'$, of $G_i$
for $i \in \brac{1,2}$. Notice that for every $r \in \brac{1,
\dots, p+1}$, the set $\brac{v_r} \cup f_1^{-1}(v_r)
\cup f_2^{-1}(v_r)$ induces a forest in $H$, which shows
$\rho(H-v) \le p+1$. Suppose now that $v \in V(G_1)$, and let
$w \in V(G_2)$ and take $G'_1 = G_1 - v$ and $G'_2 = G_2 - w$.
Let $f'_i \colon V(G'_i) \to \brac{1, \dots, p}$ be a partition
of $V(G'_i)$ into $p$ forests for $i \in \brac{1,2}$. (Such
partitions exist due to the minimality of $G_1$ and $G_2$.)
Let $f \colon V(H-v) \to \brac{1, \dots, p+1}$ be given by

\[ f(x) =
\left\{
	\begin{array}{ll}
		f_i(x)  & \mbox{if } x \in V(G'_i) \mbox{ for } i \in \brac{1,2}\\
		p+1     & \mbox{if } x \in S\cup \brac{w}
	\end{array}
\right.
\]
It is easy to see that $f$ induces a partition of $V(H)$ into
$p+1$ forests, which shows $b$.

Part $c$ follows directly from the construction of $H$.
\end{proof}

We did not succeed to obtain an analog to Corollary
\ref{joinVertex} in terms of union vertices, which would
yield an upper bound on the size of a cograph minimal
obstruction for arboricity $p$.  Instead, we derive, from
the algorithm in the next section, the following result.

\begin{theorem}
\label{thm:upbound}
	Each minimal cograph obstruction for arboricity $p$ has at
	most $O((2p)!^2)$ vertices.
\end{theorem}


\section{A polynomial algorithm}
\label{sec:alg}

The following simple observation describes the recursive
structure of $(p,q,r)$-partitions in cographs.   The only
thing to remember is that for the second statement (join
of two cographs), a forest cannot intersect both sides in
more than one vertex.

\begin{prop}
\label{rec-ob}
	\begin{enumerate}
		\item Let $G = G_u + G_d$ be a cograph with non-empty
		      subgraphs $G_u$ and $G_d$. Then $G$ has a
		      $(p,q,r)$-partition if and only if there exist
		      integers $r_u,r_d \geq 0$ such that $c=r_u+r_d$ and
		      $G_u$ and $G_d$ have a $(p,q,r_u)$-partition and a
		      $(p,q,r_d)$-partition, respectively.

		\item Let $G = G_u \oplus G_d$ be a cograph with non-empty
		      subgraphs $G_u$ and $G_d$. Then $G$ has a
		      $(p,q,r)$-partition if and only if there exist
		      integers $p_u,p_d,q_u,q_d,r_u,r_d, t_u,t_d \geq 0$
		      such that $p=p_u+p_d+ t_u+t_d$, $q=q_u+q_d$,
		      $r=r_u+r_d$ and $G_u$ and $G_d$ have a
		      $(p_u,q_u+t_d,r_u+t_u)$-partition and
		      $(p_d,q_d+t_u,r_d+t_d)$-partition, respectively.
	\end{enumerate}
\end{prop}

Notice that the integers $p_u, p_d, q_u, q_d, r_u, r_d,
t_u, t_d$ in Proposition \ref{rec-ob} are not necessarily
unique.   If $G$ is a cograph, then there exists cographs
$G_u$ and $G_d$ such that either $G = G_u + G_d$ or $G = G_u
\oplus G_d$. Suppose $G_u$ and $G_d$ have an $(x,y,z)$-partition
and an $(x',y',z')$-partition, respectively. Then, the triples
$(p,q,r)$ such that $G$ has a $(p,q,r)$-partition as the
one described in Proposition \ref{rec-ob} are said to be
{\em derived} from $(x,y,z)$ and $(x',y',z')$.

This structure can be used to obtain an efficient algorithm
to solve the $(p,q,r)$-partition problem in cographs. To this
end, define the \emph{weight} of a triple $(p,q,r)$ to be
$p+q+r$.

\begin{prop}
	\label{der-ob}
	Given two triples $T_1$ and $T_2$ with weights at most
	$m$, the set of all triples derived from $T_1$ and $T_2$
	can be generated in $O(m)$ time.
\end{prop}

\begin{proof}
If $G = G_u + G_d$ then according to the first item in
Proposition \ref{rec-ob}, we have that $(p,q,r) = (\max
\{x, x'\}, \max\{y, y'\}, z + z')$ as the only option.
But if $G = G_u \oplus G_d$ then we may have more options
for $(p,q,r)$. In this case, by setting $t = t_d+t_u$, any
triple $(x+x'+t, y+y'-t, z+z'-t)$ can be produced, given
that all the components of the triple are non-negative.
\end{proof}

\begin{theorem}
	Given a cograph $G$ with $n$ vertices, there exists
	an algorithm that computes in $O(np^7)$ all the triples
	$T$ with weight at most $p$ such that $G$ admits a
	$T$-partition.
\end{theorem}
\begin{proof}
	We build an algorithm $ALG$ recursively as follows.
	The algorithm is trivial when $G$ is a clique or an
	independent set. So suppose either $G=G_u + G_d$ or
	$G=G_u \oplus G_d$. Then apply $ALG$ to $G_1$ and
	$G_2$ to obtain the lists $L_u$ and $L_d$. Then for
	each pair $(T_u,T_d) \in L_u\times L_d$, add all the
	triples derived from $T_u$ and $T_d$ to some list $L$,
	which is our final answer. Let $f(n)$ be the run time
	of this algorithm on a cograph with $n$ vertices.
	Suppose $G,G_u$ and $G_d$ have $n$,$s$ and $n-s$
	vertices, respectively. Considering Observation
	\ref{der-ob} and the fact that lists $L_u$ and $L_d$
	each have at most $O(p^3)$ members, we get the following
	recursion:
	\[
	f(n)= f(s) + f(n-s) + O(p^7)
	\]
	which implies $f(n)=O(np^7)$.
\end{proof}


Finally, we use a similar approach, based on triples, to
prove Theorem \ref{thm:upbound}. Let $\mathcal{T}$ be a set
of triples of non-negative integers. We say a graph $G$ is a
\emph{minimal cograph obstruction for the set $\mathcal{T}$}
if the following conditions hold:
\begin{enumerate}
    \item $G$ does not admit a $T$-partition for any $T \in \mathcal{T}$.

    \item For any vertex $v \in V(G)$, there exists a triple $T
          \in \mathcal{T}$ such that $G-\{v\}$ admits a $T$-partition.
\end{enumerate}

Given a triple $T = (p, q, r)$, by the  {\em weight $w(T)$ of
$T$} now we mean $2p+q+r$ (and not $p+q+r$).   The following
observation is immediate.

\begin{obs}
	\label{weight-ob}
	Suppose a triple $T$ is derived from triples $T_u$
	and $T_d$. Then $w(T_u),w(T_d) \leq w(T)$. Furthermore,
	if $G = G_u + G_d$, $G$ is $T$-partitionable, and $G_u$
	and $G_d$ are $T_u$-partitionable and $T_d$-partitionable,
	respectively, then $w(T) = w(T_u) + w(T_d)$.
\end{obs}

For integers $k,m$, let $f(k,m)$ be the smallest integer with
the following property: any minimal obstruction with respect to
a set of triples with weight at most $k$ and at most $m$ triples
with weight exactly $k$ has size at most $f(k,m)$. Now we embark
on estimating $f(k,m)$ using recursion. Note that $m=O(k^2)$.
For the sake of convenience, we assume $f(k,0)=f(k-1,O(k^2))$.

\begin{theorem}
	\label{MO-the}
	$f(k, m) = O(k!^2)$
\end{theorem}
\begin{proof}
	Let $\mathcal{T}$ be a set of triples with weight at most $k$ and at
	most $m$ triples having weight $k$. Let $G$ a minimal
	obstruction with respect to $\mathcal{T}$. Ignoring trivial cases,
	we may assume that either $G = G_u + G_d$ or $G = G_u \oplus
	G_d$, where both $G_u$ and $G_d$ are non-empty. For $i \in
	\{ u,d \}$ denote by $L_i$ the set of triples $X$ with weight
	at most $k$ such that $G_i$ admits an $X$-partition. We say
	a triple $X$ is \emph{dangerous} for $G_u$ ($G_d$, respectively)
	if there exists a triple $X' \in L_d$ ($X' \in L_u$,
	respectively) such that from $X$ and $X'$ we can derive a
	triple in $\mathcal{T}$. Let $D_u$ ($D_d$, respectively) be the set
	of all triples dangerous for $G_u$. Note that $D_u$ is
	non-empty if $G_u$ has at least two vertices. To see this,
	let $v \in V(G_u)$ be an arbitrary vertex. Then $G - \{ v \}$
	has a $T$-partition for some $T \in \mathcal{T}$. Applying Observation
	\ref{rec-ob}, there must be triples $X$ and $X'$ such that
	$T$ is derived from $X$ and $X'$ and $G_u - \{ v \}$ and
	$G_d$ admit an $X$-partition and an $X'$-partition, respectively.
	This means $X' \in L_d$ and $X \in D_u$. A similar argument
	shows that $G_u$ must be a minimal obstruction with respect
	to $D_u$ if it has at least two vertices. Note that since $G_d$
	is non-empty so any triple in $L_d$ has non zero weight. This
	implies the weight of triples in $D_u$ is at most $k$ according
	to Observation \ref{weight-ob}. In fact, if $G = G_u \oplus G_d$,
	the weight of the triples in $D_u$ (and $D_d$) are at most $k-1$.
	So in this case, we have $|V(G_u)|,|V(G_d)| \le f(k-1,O(k^2))$,
	and thus:
	\begin{equation}
	\label{eq1}
	f(k,m) \leq 2f(k-1,O(k^2)).
	\end{equation}

	Now suppose $G=G_u + G_d$. Suppose $D_u$ contains a triple
	$X=(x, y, z)$ with weight $k$. This implies $G_d$ admits an $(x,
	y, 0)$-partition. So $X \in \mathcal{T}$, which means $X \notin D_d$
	(otherwise $G$ would admit an $X$-partition). This means that
	for some integer $t$, $0 \le t \le m$ we have
	\begin{equation}
	\label{eq2}
	f(k,m) \leq f(k,m-t)+f(k,t).
	\end{equation}

	Now (\ref{eq1}) and (\ref{eq2}) imply that $f(k,m)=O(k!^2)$.

\end{proof}



\section{Concluding remarks}
\label{sec:conc}

We have already observed that the $(p,q,r)$-partition
problem can be considered as a general framework that
includes interesting problems, e.g., $q$-colouring,
arboricity $p$, or independent feedback vertex set.
In these cases, the value of $r$ is $0$. Notice that
$r$ can be used as an additional input value to state
some classic decision problems arising from optimization
problems.   We discuss two examples.

Recall that the vertex cover optimization problem asks,
given a graph $G$, to find the size of a minimum vertex
cover of $G$.   There is a decision problem associated
to this optimization problem.   The problem
\textsc{Vertex Cover} takes as input a graph $G$ and a
non-negative integer $r$, and asks whether $G$ contains
a vertex cover with at most $r$ vertices.   Now, notice
that a $(0,1,r)$-partition of a graph $G$ is a partition
into an independent set, and a set with at most $r$
vertices $C$, this is, all the edges of $G$ must have at
least one end in the set $C$.   From here, it is easy
to conclude that $G$ has a vertex cover with at most $r$
vertices if and only if $G$ admits a $(0,1,r)$-partition.

The odd cycle transversal problem asks to find the
minimum set of vertices having a non-empty intersection
with every odd cycle in a graph $G$.   Again, this
optimization problem has an associated decision problem.
Consider the \textsc{Bipartization} problem, with input
$(G, r)$, where $G$ is a graph and $r$ is a non-negative
integer, and where we have to decide whether or not there
is a subset $X$ of at most $r$ vertices of $G$ such that
$G - X$ is a bipartite graph.   Notice that alternatively,
we could ask whether $G$ admits a $(0, 2, r)$-partition.
It was proved in \cite{lewisJCSS20} that \textsc{Bipartization}
is NP-complete, even when restricted to planar graphs.

In this setting, it is easy to notice that
$(1,0,r)$-partition corresponds to the Feedback Vertex
Set problem.   We think that $(p,q,r)$-partitions represent
a nice framework where many seemingly unrelated problems
converge.

Before proceeding, we make a simple observation about the
structure of minimal obstructions to the $(0,q,r)$-partition
problem.

\begin{theorem}
Let $q$ and $r$ be non-negative integers. Every
disconnected cograph minimal obstruction is of the form
$\bigcup_{i \in I} G_{i}$ where $G_i$ is a cograph minimal
obstruction for $(0, q, r_i)$ and $\card{I} - 1 +
\sum_{i\in I}r_i = r$.
\end{theorem}

\begin{proof}
Let $G$ be a minimal obstruction and $\brac{G_i}_{i\in I}$ the set of components of $G$. Since $G$ is a minimal obstruction, we know that for every $x \in V(G)$ there exists $L_x \subseteq V(G)$ such that $x \in L_x$, $\card{L_x} \leq r +1$, and $G-L_x$ is $q$-colourable.

\begin{claim-a}\label{Ibound}
	For every $x,y \in V(G)$, every choice of $L_x$ and $L_y$, we have $\card{L_x \cap V(G_i)} = \card{L_y \cap V(G_i)}$ for each $i \in I$.
\end{claim-a}

Suppose otherwise and let $x, y \in V(G)$, $L_x, L_y \subseteq V(G)$ and $i \in I$ be such that $\card{L_y \cap V(G_i)} < \card{L_x \cap V(G_i)}$. Since $G$ is a minimal obstruction for $(0,q,r)$-partition, then $\card{L_v}=r+1$ for every $v \in V(G)$. This means that $L_x' = \left( L_x - V(G_i)\right) \cup \left(L_y \cap V(G_i)\right)$ satisfies $\card{L_x'} \leq r$, and so $G - L_x'$ is not $q$-colourable. Since $G_j -L_x' = G_j - L_x$ for every $j \in I$, $j \neq i$, it follows that $\chi(G_i-L_x') > q$, but $G_i-L_x' = G_i-L_y$, contradicting the choice of $L_y$.  \bs

\begin{claim-a}
	Let $x \in V(G)$ and take $r_i = \card{L_x \cap V(G_i)} -1$. For every $i \in I$, the cograph $G_i$ is a minimal obstruction for $(0,q,r_i)$-partition.
\end{claim-a}

Due to the choice of $r_i$ and Claim \ref{Ibound}, $G_i$ is an obstruction to $(0,q,r_i)$-partition. To see that it is minimal, let $x \in V(G_i)$ and $L_x \subseteq V(G)$. Since $G-L_x$ is $q$-colourable and $\card{\left(V(G_i)\cap L_x\right) - x} = r_i$, then $G$ is a minimal obstruction for $(0,q,r_i)$-partition.\bs

Since $\sum_{i\in I}r_i = \card{L_x} - \card{I}$, we get that $\card{I} - 1 + \sum_{i\in I}r_i = \card{L_x}-1 = r$, completing the proof of the theorem.
\end{proof}

We propose a question that might turn out to be
interesting. Observe that for the two cases $(0,q,0)$
and $(1,q,0)$ the number of cograph minimal obstructions
is independent of $q$. (There is precisely one cograph
minimal obstruction for $(0,q,0)$-partition, because
cographs are perfect, and it follows from Theorem
\ref{ifvs} that there are exactly two cograph minimal
obstructions for $(1,q,0)$-partition.) We wonder
whether it is always the case that the number of
cograph minimal $(p,q,0)$-obstructions is independent
of $q$.

\begin{problem}
\label{probq}
Let $p$ be a fixed non-negative integer.   Is
it true that there is an integer $f(p)$ such that
for every non-negative integer $q$, the number of
cograph minimal obstructions for $(p,q,0)$-partition is
$f(p)$?
\end{problem}

If the answer to Problem \ref{probq} is negative, then
it could be interesting to determine for which values
of $p$ the number of cograph minimal obstructions
for $(p,q,0)$-partition is independent of $q$.

%
%
%
%
%
%

\end{document}